\documentclass[12pt,a4paper,reqno]{article}
\usepackage{amsmath,amsthm,amssymb}

\theoremstyle{plain} 
\newtheorem{theorem}{Theorem}
\newtheorem{proposition}{Proposition}[section]

\newtheorem{lemma}{Lemma}[section]

\theoremstyle{definition}
\newtheorem{definition}{Definition}[section]

\newtheorem*{remark}{Remark}

\theoremstyle{remark}
\newtheorem*{note}{Note}

\numberwithin{equation}{section}

\newcommand{\thmref}[1]{Theorem~\ref{#1}}
\newcommand{\secref}[1]{Section~\ref{#1}}
\newcommand{\lemref}[1]{Lemma~\ref{#1}}
\newcommand{\propref}[1]{Proposition~\ref{#1}}
\newcommand{\defref}[1]{Definition~\ref{#1}}

\newcommand{\dem}[1]{{\em{#1}}}  

\newcommand{\beq}{\begin{equation}}
\newcommand{\eeq}{\end{equation}}
\newcommand{\ot}{\otimes}
\newcommand{\hot}{\hat{\ot}\,}
\renewcommand{\t}{\tilde}
\newcommand{\id}{\text{id}}
\newcommand{\modh}{\ (\text{mod }\hhh)}
\newcommand{\ct}{\cite}
\newcommand{\half}{\frac{1}{2}}

\renewcommand{\a}{\alpha}
\renewcommand{\b}{\beta}

\renewcommand{\d}{\delta}

\newcommand{\lie}{\mathfrak{g}}	
\newcommand{\qlie}{\lie_\hhh}	
\newcommand{\lqlie}{\mathfrak{L}_\hhh(\lie)}
\newcommand{\wqlie}{\mathfrak{l}_\hhh(\lie)}   
\newcommand{\uqg}{U_\hhh(\lie)}

\newcommand{\roots}{R}
\newcommand{\sln}{\mathfrak{sl}_{n}}
\newcommand{\slnn}{\mathfrak{sl}_{n>2}}
\newcommand{\qsln}{(\sln)_\hhh}
\newcommand{\lqsln}{\mathfrak{L}_\hhh(\sln)}
\newcommand{\gln}{\mathfrak{gl}_{n}}
\newcommand{\CC}{{\mathbb C}}
\newcommand{\ch}{{\mathbb C}[[\hhh]]}
\newcommand{\pa}{\chi}   
\newcommand{\qconj}[1]{#1^\triangledown}	
\newcommand{\hhh}{h}

\newcommand{\xp}{X^+_h}
\newcommand{\xm}{X^-_h}
\newcommand{\xpm}{X^\pm_h}
\newcommand{\hh}{H_h}

\newcommand{\ads}{\text{ad}}	
\newcommand{\ad}[2]{(\ads\,{#1})\,#2}	

\newcommand{\adsh}{\text{ad}^{(\hhh)}}	
\newcommand{\adh}[2]{(\adsh\,{#1})\,#2}	
\newcommand{\adth}[2]{(\adsh_2\,{#1})\,#2}	

\newcommand{\adsc}{\text{ad}^{(0)}}	
\newcommand{\adc}[2]{(\adsc\,{#1})\,#2}	
\newcommand{\adtc}[2]{(\adsc_2\,{#1})\,#2}	

\newcommand{\qlb}{[\,,]_\hhh}	
\newcommand{\lb}{[\,,]}		
\newcommand{\driso}{\varphi}	
\newcommand{\lieh}{\lie[[\hhh]]}

\begin{document}

\title{
\vspace{-15mm}
\begin{flushright}\small q-alg/9605025 \\
\small KCL-TH-96-05\\[30pt]\end{flushright}
\Huge Quantum Lie algebras,\\
\Large their existence, uniqueness and $q$-antisymmetry}

\author{Gustav W. Delius\\[8pt] 
  \small Department of Mathematics, King's College London\\[-4pt]
  \small Strand, London WC2R 2LS, Great Britain\\[-4pt]
  \small e-mail: delius@mth.kcl.ac.uk\\[-4pt]
  \small http://www.mth.kcl.ac.uk/$\sim$delius\\[10pt]
  Mark D. Gould\\[8pt]
  \small Department of Mathematics, University of Queensland\\[-4pt]
  \small Brisbane Qld 4072, Australia}

\date{}

\maketitle

\begin{abstract}
Quantum Lie algebras are generalizations of Lie algebras which 
have the quantum parameter $h$ built into their structure.
They have been defined concretely as certain submodules $\lqlie$
of the quantized enveloping algebras $\uqg$. On them the quantum
Lie product is given by the quantum adjoint action.

Here we define for any finite-dimensional
simple complex Lie algebra $\lie$ an abstract
quantum Lie algebra $\qlie$ independent of any concrete realization.
Its $h$-dependent structure constants are given in terms of 
inverse quantum Clebsch-Gordan coefficients. 
We then show that all concrete quantum Lie algebras $\lqlie$ are
isomorphic to an abstract quantum Lie algebra $\qlie$.

In this way we prove two important properties of quantum Lie algebras:
1) all quantum Lie algebras $\lqlie$ associated to the same $\lie$
are isomorphic, 2) the quantum Lie product of any $\lqlie$ is
$q$-antisymmetric. We also describe a construction of $\lqlie$ which
establishes their existence.

\end{abstract}

\newpage

\section{Introduction}

Lie algebras play an important role in the description of many classical
physical theories. This is particularly pronounced in integrable
models which are described entirely in terms of Lie algebraic data.
However, when quantizing a classical theory the Lie algebraic
description seems to be destroyed by quantum corrections.

It is conceivable that in some cases the Lie algebraic structure of the 
theory is deformed rather than destroyed. The quantum theory may be 
describable by a quantum generalization of a Lie algebra which has 
higher order terms in $\hbar$ built into its structure. These
speculations were prompted by the beautiful structure found in
affine Toda quantum field theories \cite{toda}. They are the
physical motivation for this work on quantum Lie algebras.

As a preliminary step towards physical applications it is necessary to identify
the natural quantum generalizations of Lie algebras and to study
their properties. Quantum generalizations $\uqg$ of the enveloping algebras
$U(\lie)$ of Lie algebras $\lie$ have been known since
the work of Drinfeld \cite{Dri85} and Jimbo \cite{Jim85} and they have
been found to play a central role in quantum integrable models. This
has lead us in \cite{qlie} to define quantum Lie algebras $\lqlie$
as certain submodules of $\uqg$, modelling the way in which ordinary
Lie algebras are naturally embedded in $U(\lie)$.

Explicit examples of quantum Lie algebras were constructed in \cite{qlie} using 
symbolic computer calculations, in particular for 
$\mathfrak{L}_\hhh(\mathfrak{sl}_3)$,
$\mathfrak{L}_\hhh(\mathfrak{sl}_4)$,
$\mathfrak{L}_\hhh(\mathfrak{sp}_4)$ and
$\mathfrak{L}_\hhh(G_2)$.
It was found empirically that in these quantum Lie algebras the
quantum Lie products satisfy an intriguing generalization of the
classical antisymmetry property. They are $q$-antisymmetric.
This can be exhibited already in the simple example of
$\mathfrak{L}_\hhh(\mathfrak{sl}_2)$. This quantum Lie algebra is
spanned by three generators $\xp,\xm$ and $\hh$ with the quantum
Lie product relations
\begin{align}\label{sl2q}
[\xp,\xm]_h&=\hh,&
[\xm,\xp]_h&=-\hh,\nonumber\\
[\hh,\xpm]_h&=\pm 2 q^{\pm 1} \xpm,&
[\xpm,\hh]_h&=\mp 2 q^{\mp 1} \xpm\nonumber\\
[\hh,\hh]_h&=2(q-q^{-1}) \hh,&
[\xpm,\xpm]_h&=0.
\end{align}
Here $q=e^h$ is the quantum parameter. Clearly for $q=1$ the above reduces
to the ordinary $\mathfrak{sl}_2$ Lie algebra. For $q\neq 1$ the
Lie product is antisymmetric if the interchange of the factors is
accompanied by $q\to q^{-1}$.

To convincingly establish that the quantum Lie algebras $\lqlie$ defined
in \cite{qlie} are the natural quantum generalizations of Lie algebras,
three questions in particular should be answered:
\begin{enumerate}
\vspace{-1mm}
\item Do the $\lqlie$ exist for all $\lie$?
\vspace{-1mm}
\item Are all $\lqlie$ associated to the same $\lie$ isomorphic?
\vspace{-1mm}
\item Do all $\lqlie$ have $q$-antisymmetric quantum Lie products?
\vspace{-1mm}
\end{enumerate}
These questions will be answered in the affirmative in this paper.

The paper is organized as follows:
\secref{sec:prelim} contains preliminaries about quantized enveloping
algebras $\uqg$ and defines the concept of $q$-conjugation. In 
\secref{sec:ab} we give a new definition of quantum Lie algebras $\qlie$
which is independent of any realization as submodules of $\uqg$.
We study the properties of the $\qlie$. In \secref{sec:qlie} we recall
the definition of the quantum Lie algebras $\lqlie$ and then show that all
$\lqlie$ are isomorphic to $\qlie$. It is in this way that we arrive
in \thmref{theorem} at the answers to 
questions 2) and 3) above. In \secref{sec:constr} we describe a
construction for quantum Lie algebras $\lqlie$ for any 
finite-dimensional simple complex
Lie algebra $\lie$, thus establishing their existence.

There are many natural questions about quantum Lie algebras which we do 
not address in this paper. These are question of representations,
of the enveloping algebras, of exponentiation to quantum groups,
of applications to physics and many more which we hope will be
addressed in the future.

We do not wish to reserve the term {\it quantum Lie algebra} only
for the particular algebras defined in this paper. Rather
we view the algebras $\qlie$ and $\lqlie$ which are defined in
Definitions \ref{def:ab} and \ref{def:qlie} in terms of $\uqg$ as 
particular examples of a more general concept
of quantum Lie algebras. What a quantum Lie algebra should be in
general is not yet known, i.e., there are not yet any satisfactory
axioms for quantum Lie algebras. Finding such an axiomatic
definition is an important problem. We hope that our study of
the quantum Lie algebras arising from $\uqg$ will help to provide the
ideas needed to formulate the axioms. In particular we expect that
the $q$-antisymmetry of the product discovered here will be an 
important ingredient.

There has been an important earlier approach to the subject of
quantum Lie algebras. It was initiated by Woronowicz in his work
on bicovariant differential calculi on quantum groups \cite{Wor89}.
He defined a quantum Lie product on the dual space to the space
of left-invariant one-forms. This has been developed further by 
several groups \cite{difcalc}. These quantum Lie algebras are 
$n^2$-dimensional where
$n$ is the dimension of the defining representation 
of $\lie$ and thus they do not have the same dimension
as the classical Lie algebra except for $\lie=\gln$. It has
never been shown how to project them onto quantum Lie algebras of
the correct dimension. Only recently Sudbery \cite{Sud95} 
has defined quantum
Lie algebras for $\lie=\sln$ which have the correct dimension $n^2-1$.
These are isomorphic to our $\qsln(0)$ (set $s=-1, t=0$ in
\propref{slnstruc}). Sch\"uler and Schm\"udgen \cite{Sch} have defined 
$n^2-1$ dimensional quantum Lie algebras for $\sln$ using left-covariant
differential calculi.
In \cite{diffcalc} we explained how our quantum Lie algebras lead to 
bicovariant differential calculi of the correct dimension.

Up to date information on quantum Lie algebras 
can be found on the World Wide Web at
http://www.mth.kcl.ac.uk/$\sim$delius/q-lie.html

\section{Preliminaries}\label{sec:prelim}

We recall the definition of quantized enveloping algebras
$\uqg$ \cite{Dri85,Jim85,Cha94} in order to fix our notation.
$\uqg$ is an algebra over $\ch $, the ring of formal power 
series in an indeterminate $\hhh$. In applications of
quantum groups in physics, the parameter $\hhh$ does not need to be
identified with Planck's constant. In general it will depend on a 
dimensionless combination
of coupling constants and Planck's constant. We use the notation
$q=e^h$.

The formal power series in $\hhh$ form only a ring, not a field. 
It is not possible to divide by an element of $\ch$ unless the power
series contains a term of order $\hhh^0$. We will have to work with
modules over this ring, rather than with vector spaces over a field as
would be more familiar to physicists like ourselves. 
However $\ch$ is a principal ideal domain and
thus many of the usual results of linear algebra continue to hold 
\cite{Cur62}.

In the physics literature on quantum groups it is quite common to 
treat $q$ not as an
indeterminate but as a complex (or real) number. It is our opinion that
in doing so, physicists loose much of the potential power of quantum
groups. Keeping $\hhh$ as an indeterminate in the formalism will, when
applied to quantum mechanical systems, lead to deeper insight.

\begin{definition}
Let $\lie$ be a finite-dimensional 
simple complex Lie algebra with symmetrizable Cartan matrix
$a_{ij}$. The \dem{quantized enveloping algebra}
$\uqg$ is the unital associative algebra over $\ch$ (completed in
the $h$-adic topology) with 
generators $x_i^+,\ x_i^-,\ h_i$, 
$1 \le i \le \text{rank}(\lie)$ and relations
\footnote{Our $x_i^{\pm}$ are related to the $X_i^\pm$ of
\cite{Cha94} by $x_i^+= q_i^{-h_i/2} X_i^+$ and $x_i^-= X_i^- q_i^{h_i/2}$
and it uses the opposite Hopf-algebra structure.}
\begin{gather}\label{uqrel}
 h_i h_j = h_j h_i,~~~~~
h_i x_j^\pm - x_j^\pm h_i = \pm a_{ij} x_j^\pm , \nonumber\\
x_i^+ x_j^- -x_j^- x_i^+ =  \delta_{ij}
\ \frac{q_i^{h_i} - q_i^{-h_i}}{q_i- q_i^{-1}},\\
\sum_{k=0}^{1-a_{ij}} (-1)^k
\left[ \begin{array}{c} 1-a_{ij} \\ k \end{array} \right]_{q_i}
(x_i^\pm)^k x_j^\pm (x_i^\pm)^{1-a_{ij}-k} = 0 \qquad i \not= j.\nonumber
\end{gather}
Here $\left[\begin{array}{c}a\\b\end{array}\right]_q$ are the
q-binomial coefficients.
We have defined $q_i = e^{d_i \hhh}$ where $d_i$ are the
coprime integers such that $d_i a_{ij}$ is a symmetric matrix.
\end{definition}
The Hopf algebra structure of $\uqg$ is given by the comultiplication
$\Delta:\uqg\to\uqg\hot\uqg$ ($\hot$ denotes the tensor product over $\ch$, 
completed in the $\hhh$-adic topology when necessary) defined by
\footnote{Interchanging $q$ and $q^{-1}$ gives an alternative Hopf
algebra structure, which is the one chosen in \cite{qlie,Cha94}.}
\begin{align}
\Delta(h_i) &= h_i \hot 1 + 1 \hot h_i, \\
\Delta(x_i^\pm) &= x_i^\pm \hot q_i^{-h_i/2} +
 q_i^{h_i/2} \hot x_i^\pm,
\end{align}
and the antipode $S$ and counit $\epsilon$ defined by
\begin{equation}\label{antipode}
S(h_i)= - h_i, ~~~
S(x_i^\pm) = - q_i^{\mp 1}\,x_i^\pm ,~~~
\epsilon(h_i) = \epsilon(x_i^\pm) = 0.
\end{equation}
$\uqg$ is quasitriangular with universal $R$-matrix $R\in\uqg\hot\uqg$.
The adjoint action of $\uqg$ on itself is given, 
using Sweedler's notation \ct{Swe69}, by
\begin{equation}\label{adjoint}
\ad{x}{y}=\sum x_{(1)}\,y\,S(x_{(2)}),~~~~~x,y\in\uqg.
\end{equation}
If the Dynkin diagram of $\lie$ has a symmetry $\tau$ which maps node
$i$ into node $\tau(i)$ then $\uqg$ has a Hopf-algebra automorphism 
defined by
$\tau(x^\pm_i)=x^\pm_{\tau(i)},\ \tau(h_i)=h_{\tau(i)}$.
Such $\tau$ are referred to as diagram automorphisms and
except for rescalings of the $x^\pm_i$ they are the only
Hopf-algebra automorphisms of $\uqg$.

\begin{proposition}[Drinfel'd \cite{Dri90}]
\label{thm:dri}
There exists an algebra isomorphism 
\newline
$\driso:\uqg\to U(\lie)[[\hhh]]$
such that $\driso\equiv\text{id}\modh$ and 
$\driso(h_i)=h_i$.
\end{proposition}

\begin{note}
This is not a Hopf-algebra isomorphism however.
\end{note}

\begin{proposition}\label{thm:rep}
By $(V^\mu,\pi^\mu)$ denote the $U(\lie)$-representation with highest 
weight $\mu$, carrier space $V^\mu$ and representation map $\pi^\mu$. 
Let $\{(V^\mu,\pi^\mu)\}_{\mu\in D_+}$
be the set of all finite-dimensional irreducible representations of $U(\lie)$.
$D_+$ is the set of dominant weights. Let $m_\lambda^{\mu\nu}$ denote 
the multiplicities in the decomposition
of tensor product representations into irreducible $U(\lie)$ representations
\begin{equation}\label{clasdec}
V^\mu\ot V^\nu=\bigoplus_{\lambda\in D_+}\,m_\lambda^{\mu\nu}\,V^\lambda.
\end{equation}
Then
\begin{enumerate}
\item $\{(V^\mu[[\hhh]],\pi^\mu\circ\driso)\}_{\mu\in D_+}$ is the
set of all indecomposable representations of $\uqg$ which are 
finite-dimensional, i.e., topologically free and of finite rank. 
Here $\driso$ is the isomorphism of \propref{thm:dri}.
\item The decomposition of $\uqg$ tensor product representations into 
indecomposable $\uqg$ representations is described by the classical
multiplicities $m_\lambda^{\mu\nu}$
\begin{equation}\label{quantdec}
V^\mu[[\hhh]]\hot V^\nu[[\hhh]]=
\bigoplus_{\lambda\in D_+}\,m_\lambda^{\mu\nu}\,V^\lambda[[\hhh]].
\end{equation}
\end{enumerate}
\end{proposition}

\begin{proof}
1. is from Drinfel'd \cite{Dri90}. It follows immediately from the 
isomorphism property of $\driso$ and from the fact that the finite
dimensional representations of $U(\lie)$ have no non-trivial
deformations.
2. the decomposition can be achieved by the same method as classically.
A careful analysis shows that working over $\ch$ does not lead to 
complications. The reason is
that all expressions appearing have a non-vanishing classical term.
\end{proof}

\begin{note}
The $\uqg$ modules $V[[\hhh]]$ are not irreducible. Their submodules
are of the form $c\,V[[\hhh]]$ with $c\in\ch$ not invertible. 
In this setting Schur's lemma takes the following form:
\begin{lemma}[Schur's lemma]\label{schur}
Let $V[[\hhh]]$ and $W[[\hhh]]$ be two finite-dimen\-sio\-nal indecomposable
$\uqg$-modules and let $f:V[[\hhh]]\rightarrow W[[\hhh]]$
be a $\uqg$-module homomorphism. Then if $f\neq 0$ then $f=c\,g$ with
$c\in\ch$ and $g$ an isomorphism. 
\end{lemma}
\end{note}

A central concept in the theory of quantum Lie algebras \cite{qlie} is
$q$-conjuga\-tion which in $\ch$ maps
$\hhh\mapsto -\hhh$, i.e. $q\mapsto q^{-1}$. 

\begin{definition}
\end{definition}
\vspace{-5mm}
\begin{enumerate}
\renewcommand{\labelenumi}{(\roman{enumi})}
\item \dem{$q$-conjugation}
\mbox{$\sim: \ch \rightarrow\ch$}, $a\mapsto\t{a}$ is the
$\CC$-linear ring
automorphism defined by $\t{\hhh}=-\hhh$.
\vspace{-1mm}
\item Let $M,N$ be $\ch$-modules. An additive map
$\phi:M\rightarrow N$ is said to be \dem{$q$-linear} if
$\phi(\lambda \,a)=\t{\lambda}\,\phi(a),~\forall a\in M, \lambda\in\ch$.
\vspace{-1mm}
\item A \dem{$q$-conjugation on a $\ch$ module $M$} is a
$q$-linear involutive map $\qconj{}:M\rightarrow M$ with
$\qconj{}=\id\modh$.
\end{enumerate}
Note the analogy between the concepts of $q$-conjugation and complex
conjugation and between $q$-linear maps and anti-linear maps.
\begin{remark}
If $M$ is a finite-dimensional $\ch$-module then a $q$-conjugation
$\qconj{}$ on $M$ is uniquely specified by giving a basis $\{b_i\}$
which is invariant. Then the $q$-conjugation takes the form
$\qconj{(\sum_i \lambda_i b_i)}=\sum\t{\lambda}_i b_i$. 
Conversely, for any $q$-conjugation on $M$
there exists an invariant basis. It can be constructed from an arbitrary
basis by adding correction terms order by order in $h$.
\end{remark}

The unique $q$-linear algebra
auto\-mor\-phism \mbox{$\sim: \uqg \rightarrow \uqg$} which 
extends $q$-conjugation on $\ch$
by acting as the identity on the generators $x_i^\pm$ and $h_i$
is a $q$-conjugation on $\uqg$.
It exists because the relations \eqref{uqrel}
are invariant under $q\mapsto q^{-1}$.
We choose the isomorphism $\driso$ in \propref{thm:dri} such that
$\sim\circ\,\driso=\driso\,\circ\sim$. This
$q$-conjugation is a coalgebra q-antiautomorphism of $\uqg$, i.e.,
$\epsilon\, \circ \sim = \sim \circ \,\epsilon,~~
\Delta\, \circ \sim = \sim \circ\, \Delta^T$ and it satisfies
$S \,\circ \sim = \sim \circ \,S^{-1}$. The map $\sim$ was introduced
already in \cite{Dri90}.

If in physical applications
$q$ were identified with a combination of a coupling constant and Planck's
constant, then $q$-conjugation would correspond to the strong-weak coupling
duality\footnote{
In some applications of quantum groups the relation between $q$ and
the coupling constant is not linear but exponential and then
$q$-conjugation is not related to strong-weak duality}. 
It has been observed in several quantum field theories, that
such a duality transformation can form a symmetry of the theory.
Affine Toda field theories in two dimensions \cite{toda} as well
as supersymmetric Yang-Mills theory in four dimensions provide 
examples of this phenomenon. It is thus very desirable to have an
algebraic structure, in which $q$-conjugation is incorporated. 
We hope that the study of this structure will one day enhance our
understanding of the origin of strong-weak coupling duality in physics.
\section{Quantum Lie algebras $\qlie$}\label{sec:ab}

The quantized enveloping algebra $\uqg$ is an infinite dimensional algebra.
It is our aim to associate to it in a natural way a finite dimensional 
algebra which would be the quantum analog of the Lie algebra.
Here our approach is based on the observation that classically a 
Lie algebra $\lie$ is also the carrier space of the adjoint
representation $\adsc$ of $U(\lie)$. The superscript $0$ is to remind
us that this is the classical adjoint representation. It is defined by
$\adc{a}{b}=[a,b]\ \forall a,b\in\lie$. It follows from the Jacobi identity
\begin{equation}
[a,[b,c]]=[[a,b],c]+[b,[a,c]]~~~ \forall a,b,c\in\lie
\end{equation}
that
\begin{equation}\label{succ}
(\adsc\,x)\circ \lb=\lb\circ(\adsc_2\,x),~~~~
\forall x\in U(\lie),
\end{equation}
where $(\adsc_2\,x)=(\adsc\,\ot \adsc)\,\Delta(x)$
is the tensor product representation carried by $\lie\ot\lie$.
Equation \eqref{succ} states that the Lie product $\lb$ of $\lie$ is 
a $U(\lie)$-module homomorphism from $\lie\ot\lie$ to $\lie$.

Because of \propref{thm:rep} we know that $\lieh$
is an indecomposable module of $\uqg$. 
Let us denote the representation of $\uqg$ on $\lieh$ by
$\adsh$. Note that at this point there is no relation between the
representation $\adsh$ of $\uqg$ on $\lieh{}$ and the adjoint
action $\ads$ of $\uqg$ on $\uqg$ defined in \eqref{adjoint}.
Generalizing the above classical observation 
we obtain a natural definition for a quantum Lie algebra
\footnote{As Ding has informed us, he and Frenkel have been 
pursuing similar ideas for some time. See also their paper
\cite{Din94} in which the utility of defining algebraic structures
using $\uqg$ module homomorphisms is stressed.}.

\begin{definition}\label{def:ab}
Let $\qlb{}:\lieh{}\hot\lieh{}\to\lieh{}$
be a $\uqg$-module homomorphism which satisfies
$\qlb=\lb\modh$. $\qlb$ gives $\lieh{}$ the structure of a non-associative
algebra over $\ch$.
We call this algebra
$\qlie=(\lie[[h]],\qlb{})$ a \dem{quantum Lie algebra}
and the product $\qlb$ a \dem{quantum Lie product}.
\end{definition}

For each Lie algebra $\lie$ this definition potentially gives many
different quantum Lie algebras $\qlie$, one for each choice
of the homomorphism $\qlb{}$.
This would be unsatisfactory were it not for the fact that such
a $\uqg$-module homomorphism is almost unique.

\begin{proposition}\label{thm:iso}
For a given $\lie\neq\slnn$ the quantum Lie algebra $\qlie$ is unique 
(up to a rescaling of the product by an invertible element of $\ch$). 
For $\lie=\sln$ with $n>2$ there is a family of quantum
Lie algebras $\qsln(\pa)$ depending on a parameter $\pa\in\CC((\hhh))$
(see \propref{slnstruc}). 
\end{proposition}

\begin{proof}
The idea of the proof is simple: For $\lie\neq\slnn$
the adjoint representation appears in the tensor product of two 
adjoint representations with unit multiplicity.
This is an empirical fact. Thus the homomorphism
$\qlb$ from $\lieh\hot\lieh$ into $\lieh$ with the requirement that
$\qlb\modh=\lb$ is unique by the weak form of Schur's lemma.

In the case $\lie=\sln$ with $n>2$ however, 
the adjoint representation appears
with multiplicity two in the tensor product. Any module arising from
a linear combination of the highest weight vectors of 
two adjoint modules is also an adjoint module and this leads to
a one-parameter family of non-isomorphic weak quantum Lie algebras 
$\qsln(\pa)$.

We find it helpful to be more explicit here than necessary and to explain
how the homomorphism $\qlb$ is obtained from inverse Clebsch-Gordan
coefficients. We begin with $\lie\neq\slnn$ and with the classical 
situation.

Let $\{v_a\}$ be a basis for $\lie$ which contains a highest weight
vector $v_0$, i.e.,
\begin{equation}\label{highest}
\adc{x^+_i}{v_0}=0,~~~~
\adc{h_i}{v_0}=\psi(h_i)v_0,~~~
\forall i,
\end{equation}
where $\psi$ is the highest root of $\lie$. Let $P_a(x^-)$ be the
polynomials in the $x^-_i$ such that $v_a=\adc{P_a(x^-)}{v_0}$.
The adjoint representation matrices $\pi$ in this basis are defined by
\begin{equation}\label{admat}
\adc{x}{v_a}=v_b\,\pi^b{}_a(x).
\end{equation}
In this paper we use the summation convention according to which
repeated indices are summed over their range.

$\lie\ot\lie$ contains a highest weight state $\hat{v}_0$ such that
\begin{equation}\label{highest2}
\adtc{x^+_i}{\hat{v}_0}=0,~~~~
\adtc{h_i}{\hat{v}_0}=\psi(h_i)\hat{v}_0,~~~
\forall i,
\end{equation}
For $\lie\neq\slnn$ this state is unique up to rescaling.
The vectors
\begin{equation}\label{basis2}
\hat{v}_a=\adtc{P_a(x^-)}{\hat{v}_0}=K_a{}^{bc}\,v_b\ot v_c
\end{equation}
form a basis for $\lie$ inside $\lie\ot\lie$ such that
\begin{equation}\label{admat2}
\adtc{x}{\hat{v}_a}=\hat{v}_b\,\pi^b{}_a(x)
\end{equation}
with the same representation matrices $\pi$ as in \eqref{admat}.
Thus the map 
\begin{equation}\label{beta}
\beta:v_a\mapsto\hat{v}_a=K_a{}^{bc}\,v_b\ot v_c
\end{equation}
is a $U(\lie)$-module homomorphism $\beta:\lie\to\lie\ot\lie$.
The coefficients $K_a{}^{bc}$ are called the Clebsch-Gordan
coefficients.
$\lie$ and Im$(\beta)$ are irreducible modules and thus by Schur's
lemma the homomorphism $\beta$ is invertible on its image.
Define $\lb:\lie\ot\lie\to\lie$ 
to be zero on the module complement of the image of $\beta$ and
on the image of $\beta$ define
$\lb=\beta^{-1}$. Then $\lb$ is the $U(\lie)$ homomorphism
from $\lie\ot\lie$ to $\lie$, unique up to rescaling. 
It is the Lie product of $\lie$.
On the basis it is given by
\begin{equation}\label{cs}
[v_a,v_b]=f_{ab}{}^c \,v_c,~~~~
\text{where }K_a{}^{bc}f_{bc}{}^d=\d_a{}^d.
\end{equation}
Thus the structure constants are given by the inverse Clebsch-Gordan
coefficients.

We turn to the quantum case. Let $\hat{v}_0$ be a highest weight state
inside $\lieh\hot\lieh$ satisfying the analog of \eqref{highest2} 
\begin{equation}\label{highest3}
\adth{x^+_i}{\hat{v}_0}=0,~~~~
\adth{h_i}{\hat{v}_0}=\psi(h_i)\hat{v}_0,~~~
\forall i,
\end{equation}
where $\adsh$ is the deformed adjoint representation $\adsh=\adsc\circ\driso$
and $\hat{v}_0\modh\neq 0$.
$\hat{v}_0$ generates the $\uqg$ module $\lieh$ inside $\lieh\hot\lieh$.
$\hat{v}_0$ must be unique up to rescaling, otherwise $\lieh$ would 
appear with multiplicity greater than one in $\lieh\hot\lieh$.
We construct a basis $\{\hat{v}_a\}$ as in \eqref{basis2} using
$P_a(x^-)\in\uqg$ with the
same polynomials $P_a$ as in \eqref{basis2}. This leads to quantum 
Clebsch-Gordan
coefficients $K_a{}^{bc}(\hhh)\in\ch$. We obtain a $\uqg$-module
homomorphism $\beta:\lieh\to\lieh\hot\lieh$ as in \eqref{beta}.

$\beta$ is invertible by the weak form of Schur's lemma. A homomorphism
$\qlb:\lieh\hot\lieh\to\lieh$ is obtained as above \eqref{cs}
\begin{equation}\label{qs}
[v_a,v_b]_\hhh=f_{ab}{}^c(\hhh) \,v_c,~~~~
\text{where }K_a{}^{bc}(\hhh)f_{bc}{}^d(\hhh)=\d_a{}^d.
\end{equation}
Up to rescaling it is the unique such homomorphism with the 
property that $\qlb\modh\neq 0$.

We now turn to $\lie=\sln$ with $n>2$ and again begin by considering
the classical situation. There are two linearly
independent highest weight vectors $\hat{v}^{(+)}_0$ and $\hat{v}^{(-)}_0$
in $\lie\ot\lie$ which satisfy \eqref{highest2}. They can be chosen so
that
\begin{equation}\label{parity}
\sigma\,\hat{v}_0^{(\pm)}=\pm\,\hat{v}_0^{(\pm)},
\end{equation}
where $\sigma$ is the bilinear map acting as $\sigma\,(v_a\ot v_b)=
v_b\ot v_a$. Expressed differently, the Clebsch-Gordan coefficients
$K^{(\pm)}_a{}^{bc}$ defined as in \eqref{basis2} satisfy
$K^{(\pm)}_a{}^{bc}=\pm K^{(\pm)}_a{}^{cb}$. Any linear combination
of $\hat{v}^{(+)}_0$ and $\hat{v}^{(-)}_0$ is a highest weight state
and leads to a homomorphism as described above but clearly only 
$\hat{v}^{(-)}_0$ leads to an {\em{antisymmetric}} Lie product.

In the quantum case too there are two linearly independent highest weight
states satisfying \eqref{highest3}. 
We can choose any linear combination and thus have
a one-parameter family of $\hat{v}_0(\pa)=K_0{}^{bc}(\pa,\hhh)\,
(v_b\hot v_c)$. We impose $\hat{v}_0(\pa)\modh\neq 0$ as before. 
In this way we obtain the family 
$\qsln(\pa)$ of quantum Lie algebras.
We will give these explicitly in \propref{slnstruc}. Certain
values for $\pa$ will lead to a $q$-antisymmetric quantum Lie product
(see \propref{thm:anti}).
\end{proof}

Some important properties of $\lie$ carry over immediately to $\qlie$.
Define root subspaces $\lie^{(\a)}$ of $\lie$ by
\begin{equation}
\lie^{(\a)}=\{ x\in\lie|\adc{h_i}{x}=\a(h_i)\,x\ \forall i\}.
\end{equation}
$\lie$ possesses a gradation
\begin{equation}
\lie=\bigoplus_{\a\in\roots\cup\{0\}}\lie^{(\a)}
,\qquad\qquad
[\lie^{(\a)},\lie^{(\b)}]\subset\lie^{(\a+\b)},
\end{equation}
where $\roots$ is the set of non-zero roots of $\lie$. 

\begin{proposition}\label{thm:grad}
A quantum Lie algebra $\qlie$ possesses a gradation
\begin{equation}\label{gradation}
\qlie=\bigoplus_{\a\in\roots\cup\{0\}}\lie^{(\a)}[[\hhh]],\qquad\qquad
\left[\lie^{(\a)}[[\hhh]],\lie^{(\b)}[[\hhh]]\right]_\hhh{}
\subset\lie^{(\a+\b)}[[\hhh]].
\end{equation}
\end{proposition}

\begin{proof}
According to \propref{thm:dri} the algebra isomorphism
$\driso:\uqg\to U(\lie)[[\hhh]]$ leaves the $h_i$ invariant and
thus
\begin{equation}
\lie^{(\a)}[[\hhh]]=\{ x\in\lieh|\adh{h_i}{x}=\a(h_i)\,x\ \forall i\}.
\end{equation}
Let $X_\a\in\lie^{(\a)}[[\hhh]]$ and $X_\b\in\lie^{(\b)}[[\hhh]]$.
From the homomorphism property of $\qlb$ and the coproduct 
$\Delta(h_i)=h_i\hot 1+1\hot h_i$ it follows that
\begin{align}
\adh{h_i}{[X_\a,X_\b]_\hhh}&=\left[\adh{h_i}{X_\a},X_\b\right]_\hhh
+\left[X_\a,\adh{h_i}{X_\b}\right]_\hhh
\nonumber\\
&=\left(\a(h_i)+\b(h_i)\right)\,[X_\a,X_\b]_\hhh
\end{align}
and thus $[X_\a,X_\b]_\hhh\in\lie^{(\a+\b)}[[\hhh]]$.
\end{proof}

Choosing basis vectors $X_\a\in\lie^{(\a)}$ and $H_i\in\lie^{(0)}$ 
\propref{thm:grad} implies that the quantum Lie product relations are of 
the form
\begin{alignat}{2}\label{qliestruc}
{[}H_i,X_\a]_\hhh{}&=l_\a(H_i)\,X_\a,~~~&
{[}X_\a,H_i]_\hhh{}&=-r_\a(H_i)\,X_\a,\nonumber\\
{[}H_i,H_j]_\hhh{}&=f_{ij}{}^k\,H_k,&
{[}X_\a,X_{-\a}]_\hhh{}&=g_\a{}^k\,H_k,\\
{[}X_\a,X_\b]_\hhh{}&=N_{\a\b}\,X_{\a+\b}
&\text{ for } \a+\b\in\roots, &~~~0\text{ otherwise}.\nonumber
\end{alignat}

This is similar in form to the Lie product relations of ordinary
Lie algebras. The most important differences are
\begin{enumerate}
\item The structure constants are now elements of $\ch$, i.e., they
depend explicitly on the quantum parameter.
\item $[H_i,H_j]_\hhh{}$ does not have to be zero. Thus the grade zero
subalgebra $\lie^{(0)}[[\hhh]]$ of $\qlie$ is not abelian. We will nevertheless
refer to it as the quantum Cartan subalgebra.
\item Each classical root $\a$ splits up
into a ``left'' root $l_\a$ and a ``right'' root $r_\a$. Classically
they are forced to be equal because of the antisymmetry of the Lie
product.
\end{enumerate}

The quantum Clebsch-Gordan coefficients which describe 
the homomorphism $\qlb{}: \qlie\hot\qlie\to\qlie$ can be
calculated directly by decomposing the tensor product representation.
This is however very tedious in general.
In \cite{qgln} it was done for $\qsln$ in an indirect way by 
using the R-matrix of $U_q(\sln)$. The method is based on realizing the
quantum Lie algebra as a particular submodule of $\uqg$ as explained
in \secref{s:concrete}. The particular submodule used in \cite{qgln} gives
the quantum Lie algebra $\qsln(\pa=1)$
but the method can be extended and gives the following result.

\begin{proposition}\label{slnstruc}
The parameter $\pa\in\CC((h))$ of $\qsln(\pa)$ is a fraction $\pa=t/s$ with
$s,t\in\ch$ and with the restriction that $(s+t)^{-1}\in\ch$.
The Lie product relations for $\qsln(\pa)$ are
\begin{gather}
\label{stb}
[H_k,X_{ij}]_\hhh{}=l_{ij}(H_k)\,X_{ij},~~~~~~
[X_{ij},H_k]_\hhh{}=-r_{ij}(H_k)\,X_{ij},\nonumber\\{}
[H_i,H_j]_\hhh{}=f_{ij}{}^k\,H_k,~~~~~~
[X_{ij},X_{ji}]_\hhh{}=g_{ij}{}^k\,H_k,\\{}
[X_{ij},X_{kl}]_\hhh{}=\d_{jk}\d_{i\neq l}N_{ijl}\,X_{il}-
\d_{il}\d_{j\neq k}M_{kij}\,X_{kj},\nonumber
\end{gather}
where $\{X_{ij}\}_{i,j=1\cdots n}\cup\{H_i\}_{i=1\cdots n-1}$ is a
basis and the structure constants are explicitly given by
\begin{align}
l_{ij}(H_k)&=
(q^{1-k}\d_{ki}-q^{-1-k}\d_{k,i-1})(s+t\, q^{n})\nonumber\\
&\quad\quad-(q^{k-1}\d_{kj}-q^{k+1}\d_{k,j-1})(s+t\,q^{-n}),
\\
r_{ij}(H_k)&=-l_{ji}(H_k),\label{lr}
\\
f_{ij}{}^k&=
\d_{ij}\,\left(\d_{ki}\left(s\,(q^{k+1}-q^{-k-1})+
t\,(q^{n+1-i}-q^{-n-1+i})\right)\right.
\nonumber\\
&\qquad\qquad \left.+s\,\d_{k<i}\,(q+q^{-1})(q^k-q^{-k})\right.
\nonumber\\
&\qquad\qquad \left.+t\,\d_{k>i}\,(q+q^{-1})(q^{n-k}-q^{-n+k})\right)
\nonumber\\
&\qquad+\d_{i,j-1}\left(s\,\d_{k\leq i}\,(q^{-k}-q^k)
+t\,\d_{k>i}\,(q^{k-n}-q^{-k+n})\right)
\nonumber\\
&\qquad+\d_{j,i-1}\left(s\,\d_{k\leq j}\,(q^{-k}-q^k)
+t\,\d_{k>j}\,(q^{k-n}-q^{-k+n})\right),
\\
g_{ij}{}^k&=q^{i-j}\left(
s\,\left(q^k\,\d_{k<j}-q^{-k}\,\d_{k<i}\right)+
t\,\left(q^{n-k}\d_{k\geq i}-q^{k-n}\d_{k\geq j}\right)
\right)
\nonumber\\
N_{ijl}&=q^{1/2-j}\left(s+t\,q^n\right),~~~~~~~
M_{kij}=q^{i-1/2}\left(s+t\,q^{-n}\right)
\end{align}
(We use a generalized Kronecker delta notation, e.g., $\d_{i\leq j}=1$
if $i\leq j$, $0$ otherwise.)
\end{proposition}

The restriction that if $\pa$ is written as $\pa=t/s$ then $s+t$ has
to be invertible comes from the requirement that the quantum Lie
product should not vanish modulo $h$.
For details of the calculation leading to the above formulae we refer
the reader to \cite{qgln}.

The Lie algebra $\sln$ with $n>2$ possesses an automorphism which is due to
the symmetry of the Dynkin diagram. It would be natural to require
that this automorphism survives also at the quantum level.
By inspecting the above Lie product relations we find

\begin{proposition}\label{thm:t}
The quantum Lie algebra $\qsln(\pa)$ possesses the Dynkin diagram
automorphism 
\begin{equation}
\tau(X_{ij})=-X_{n+1-j,n+1-i},~~~~~
\tau(H_i)=H_{n-i}
\end{equation}
iff $\pa=1$.
\end{proposition}
This is the reason why in \cite{qgln} we focused our attention on
the case of $\pa=1$.

The most basic property of a Lie product is its antisymmetry. In
quantum Lie algebras this has found an interesting generalization.

\begin{proposition}\label{thm:anti}
The quantum Lie product of $\qlie$ for $\lie\neq\slnn$ and of $\qsln(\pa)$ 
with $\tilde{\pa}=\pa$ is \dem{$q$-antisymmetric}, i.e.,
there exists a q-conjugation $\qconj{}:\qlie\to\qlie$ consistent with the
gradation \eqref{gradation} such that
\begin{equation}
\qconj{[a,b]_\hhh}{}=-[\qconj{b},\qconj{a}]_\hhh{}.
\end{equation}
Thus, choosing the basis in \eqref{qliestruc} so that
$\qconj{X_\a}=X_\a$, $\qconj{H_i}=H_i$, the structure constants satisfy
\begin{equation}
r_\a=\t{l}_\a,~~~f_{ij}{}^k=-\t{f}_{ji}{}^k,~~~
g_\a{}^k=-\t{g}_{-\a}{}^k,~~~N_{\a\b}=-\t{N}_{\b\a}.
\end{equation}
\end{proposition}
\begin{proof}
For $\qsln$ the statement can be verified directly from the expressions
in \propref{slnstruc}. For $\lie\neq\sln$ we use the same notation as
in the proof of \propref{thm:iso}. The adjoint representation
appears with multiplicity one in the tensor product and thus we know that
the highest weight state $\hat{v}_0=K_0{}^{ab}(\hhh)\,v_a\ot v_b$ in
$\lieh\hot\lieh$ satisfying \eqref{highest3} is unique up to rescaling.
$\t{\hat{v}}_0^T=K_0{}^{ba}(-\hhh)\,v_a\ot v_b$ also satisfies the highest
weight condition \eqref{highest3}.
\begin{align}
\adth{x_i^+}{\t{\hat{v}}_0^T}&=
\left((\adsh\ot\adsh)\Delta(x_i^+)\right)\,\t{\hat{v}}_0^T
\nonumber\\
&=\,\sim\left[\left((\adsh\ot\adsh)\Delta^T(x_i^+)\right)\,
\hat{v}_0^T\right]
\nonumber\\
&=\,\sim\left[\adth{x_i^+}{\hat{v}_0}\right]^T
\nonumber\\
&=0.
\end{align}
We used that $\sim\circ\,(\adsh\,x)=(\adsh\,\t{x})\,\circ\sim$
(which follows from $\sim\circ\,\driso=\driso\,\circ\sim$),
that $\t{v}_a=v_a$ and that $\sim\circ\,\Delta=\Delta^T\,\circ\sim$.
Thus $\hat{v}^\prime_0=\half(\hat{v}_0-\t{\hat{v}}_0^T)$ is a
highest weight state (proportional to $\hat{v}_0$ by uniqueness). 
It is non-zero because it is non-zero classically.
Following a similar calculation to the above one finds that it leads
to Clebsch-Gordan coefficients 
$K^\prime_a{}^{bc}(\hhh)=\half\left(K_a{}^{bc}(\hhh)-
K_a{}^{cb}(-\hhh)\right)$.
These are manifestly $q$-antisymmetric. Following through the
construction of the structure constants one finds
$f^\prime_{ab}{}^c(\hhh)=-f^\prime_{ba}{}^c(-\hhh)$.
\end{proof}

\section{Quantum Lie algebras $\lqlie$ inside $\uqg$}\label{s:concrete}
\label{sec:qlie}

In \defref{def:ab} quantum Lie algebras are defined abstractly, i.e.,
independently of any specific realization. 
In \cite{qlie} quantum Lie algebras were defined
as concrete objects, namely as certain submodules of the quantized
enveloping algebras $\uqg$. This definition is based on the
observation that an ordinary Lie algebra $\lie$ can be naturally
viewed as a subspace of its enveloping algebra $U(\lie)$ with the Lie
product on this subspace given by the adjoint action of $U(\lie)$.
Thus it is natural to define a quantum Lie algebra as an analogous
submodule of the quantized enveloping algebra $\uqg$ with the quantum
Lie product given by the adjoint action of $\uqg$. Before we can
state the precise definition we need some preliminaries.

The Cartan involution $\theta:\uqg\to\uqg$ 
is given by the same formulas as in the
classical case:
$\theta(x_i^\pm)=x_i^\mp,\ \theta(h_i)= -h_i$.
It is an algebra automorphism and a coalgebra antiautomorphism, i.e.,
$\Delta\circ\theta=(\theta\hot\theta)\circ\Delta^T$ and
$S\circ \theta=\theta\circ S^{-1}$.
We define a tilded Cartan involution by composing the Cartan 
involution with $q$-conjugation, i.e., 
$\t{\theta}=\sim\circ\theta$. Similarly we define a tilded
antipode as $\t{S}=\sim\circ S$.
With respect to the adjoint action defined in \eqref{adjoint} they
satisfy
$\ad{\t{\theta}(a)}{\t{\theta}(b)}=\t{\theta}(\ad{a}{b})$ and
$\ad{\t{S}(a)}{\t{S}(b)}=\t{S}(\ad{S^{-1}(a)}{b})$ for all $a,b\in\uqg$.

\begin{definition}\label{def:qlie}
A \dem{quantum Lie algebra $\lqlie$ inside $\uqg$}
is a finite-dimen\-sio\-nal 
indecomposable $\text{ad}$ - submodule
of $\uqg$ endowed with the 
\dem{quantum Lie product} 
$[a,b]_\hhh{}=\ad{a}{b}$
such that
\begin{enumerate}
\vspace{-2mm}
\item $\lqlie$ is a deformation of $\lie$, i.e., there is an
algebra isomorphism
$\lqlie\cong\lie\modh$.
\vspace{-2mm}
\item
$\lqlie$ is invariant under $\t{\theta}$, $\t{S}$ and any diagram
automorphism $\tau$.
\end{enumerate}
\vspace{-2mm}
A \dem{weak quantum Lie algebra $\wqlie$} is defined similarly but
without the requirement 2.
\end{definition}
The existence of a Cartan involution and
an antipode on $\lqlie$ plays an important role in the investigations into
the general structure of quantum Lie algebras in \cite{qlie}. In particular
it allows the definition of a quantum Killing form.
The invariance under the diagram automorphisms $\tau$ is less important
but is clearly a natural condition to impose.
It is shown in \cite{qlie} that given any weak quantum Lie algebra $\wqlie$
inside $\uqg$, one can always construct a true quantum Lie
algebra $\lqlie$ which satisfies property 2 as well. Thus this
extra requirement is not too strong. 

We now come to the relation between the abstract quantum Lie algebras
$\qlie$ of \defref{def:ab} and the concrete weak quantum Lie algebras
$\wqlie$ of \defref{def:qlie}.

\begin{proposition}\label{thm:lq}
All weak quantum Lie algebras
$\wqlie$ inside $\uqg$ are isomorphic to the quantum Lie algebra
$\qlie$ as algebras
(or to $\qsln(\pa)$ for
some $\pa$ in the case of $\lie=\sln$).
\end{proposition}

\begin{proof}
By definition $\wqlie$ is a finite-dimensional, indecomposable $\uqg$ 
module. Condition
1 of the definition implies that the representation of $\uqg$ carried by
this module is a deformation of the representation of $U(\lie)$
carried by $\lie$. There is only one such deformation, namely the
adjoint representation $\adsh$ carried by $\lieh$. Thus $\wqlie$
is isomorphic to $\lieh{}$ as a $\uqg$ module.
The identity
\begin{equation}
\sum\ad{\ad{x_{(1)}}{a}}{(\ad{x_{(2)}}{b})}=\ad{x}{(\ad{a}{b})}
\end{equation}
can be rewritten using that, when restricted to $\wqlie\subset\uqg$, 
$[a,b]_h=\ad{a}{b}=\adh{a}{b}$.
\begin{equation}
(\adsh\,x)\circ \qlb=\qlb\circ(\adsh_2\,x),~~~~
\forall x\in U(\lie).
\end{equation}
This states that the quantum Lie product on $\wqlie$ is 
a $\uqg$-module homomorphism
and thus is a quantum Lie product in the sense of \defref{def:ab}.
\end{proof}

\begin{remark} One should not confuse the adjoint {\em{action}} $\ads$
with the adjoint {\em{representation}} $\adsh$. The adjoint action
$\ads$ is defined using the coproduct and the antipode as
\[
\ad{x}{y}=x_{(1)}yS(x_{(2)})~~~~\forall x,y\in\uqg.
\] 
The adjoint
representation $\adsh$ is defined using the algebra isomorphism
$\driso:\uqg\to U(\lie)[[\hhh]]$ of \propref{thm:dri} as
\[
\adh{x}{a}=\adc{\driso(x)}{a}~~~~\forall
x\in\uqg,\ a\in\lieh.
\]
Thus the adjoint action is determined by the $h$-deformed Hopf-algebra
structure whereas the adjoint representation is determined by only
the $h$-deformed algebra structure. From this point of view it is
surprising that the two ever coincide. But the weak quantum Lie
algebras $\wqlie$ are exactly those embeddings of $\lieh$ into
$\uqg$ on which $\ads$ and $\adsh$ coincide and we will establish
their existence in the next section.
\end{remark}

\propref{thm:lq} allows us 
to answer two important questions about the
concrete quantum Lie algebras $\lqlie$ inside $\uqg$ which were left
unanswered in \cite{qlie}.

\begin{theorem}\label{theorem}
Given any finite-dimensional simple complex Lie algebra $\lie$.
\begin{enumerate}
\item
All quantum Lie
algebras $\lqlie$ are isomorphic as algebras.
\item
All quantum Lie algebras $\lqlie$ have $q$-antisymmetric Lie products.
\end{enumerate}
\end{theorem}

\begin{proof}
1. For $\lie\neq\slnn$ this is obvious from \propref{thm:lq} and the
uniqueness of $\qlie$ according to \propref{thm:iso}. For $\lie=\slnn$
the requirement of $\tau$-invariance in \defref{def:qlie} implies through
\propref{thm:t} that $\lqsln$ can be isomorphic only to $\qsln(\pa=1)$.
2. This is obvious because $\qlie$ and $\qsln(\pa=1)$ have
$q$-antisymmetric Lie products according to \propref{thm:anti}.
\end{proof}

\section{Construction of quantum Lie algebras $\lqlie$}\label{sec:constr}

There is a general method for the construction of weak quantum Lie
algebras $\wqlie$ and quantum Lie algebras $\lqlie$
inside $\uqg$. The method was presented in \cite{qgln}
for $\lie=\sln$ but it works for any finite-dimensional 
simple complex Lie algebra $\lie$ as we will discuss here.

We begin with a lemma giving a construction of ad-submodules of
$\uqg$.

\begin{lemma}\label{lem}
Let $A$ be any element of $\uqg\hot\uqg$ satisfying 
$A\,\Delta(x)=\Delta(x)\,A,~\forall x\in\uqg$. Let $V[[h]]$
be any finite-dimensional indecomposable $\uqg$ module 
and let $\pi_{ij}$ be the corresponding
representation matrices. Then the elements
\begin{equation}
A_{ij}=(\pi_{ij}\ot\id)\,A\in\uqg
\end{equation}
span an ad-submodule of $\uqg$ which is isomorphic to a
submodule of\hfill\newline
$V[[h]]^*\hot V[[h]]$, i.e.,
\begin{equation}
\ad{x}{A_{ij}}=A_{kl}\,\pi_{ki}^*(x_{(1)})\,\pi_{lj}(x_{(2)}),~~~
\forall x\in\uqg.
\end{equation}
Here $\pi^*$ denotes the dual (contragredient) representation to $\pi$
defined by
\begin{equation}\label{dualrep}
\pi_{ki}^*(x)=\pi_{ik}(S(x)). 
\end{equation}
\end{lemma}

\begin{proof}
We first calculate
\begin{align}
x\,A_{ij}&=\left(\pi_{ij}\ot\id\right)(1\ot x)\ A
\nonumber\\
&=\left(\pi_{ij}\ot\id\right)(S(x_{(1)})\ot 1)
\ A\ (x_{(2)}\ot x_{(3)})
\\
&=\pi_{ik}(S(x_{(1)}))\,A_{kl}\,\pi_{lj}(x_{(2)})\,x_{(3)}.
\nonumber
\end{align}
Then, using \eqref{dualrep}
\begin{align}
\ad{x}{A_{ij}}&=x_{(1)}\,A_{ij}\,S(x_{(2)})
\nonumber\\
&=A_{kl}\,\pi^*_{ki}(x_{(1)})\,\pi_{lj}(x_{(2)})\,x_{(3)}\,
S(x_{(4)})
\\
&=A_{kl}\,\pi^*_{ki}(x_{(1)})\,\pi_{lj}(x_{(2)})
\nonumber
\end{align}
\end{proof}

This lemma can be applied to construct weak quantum Lie algebras.

\begin{proposition}
Let $A=\hhh^{-1}\left(R^T R-1\right)$ where $R$ is the universal
R-matrix of $\uqg$ and $R^T$ the same with the tensor factors
interchanged (i.e., if $R=\sum a_i\ot b_i$ then 
$R^T=\sum b_i\ot a_i$).
Let $\{e_i\}$ be a basis for the $\uqg$ module $V[[h]]$ 
and let $\pi_{ij}$ be the corresponding representation
matrices. Choose a basis $\{v_a\}$ for the adjoint representation
$\lie[[h]]$ of $\uqg$ and let
$K:\lie[[h]]\to V[[h]]^*\hot V[[h]],
v_a\mapsto\hat{v}_a=K_a{}^{ij}\,(e^*_i\ot e_j)$ be a
$\uqg$-module homomorphism, i.e., the $K_a{}^{ij}$ are 
quantum Clebsch-Gordan coefficients.
Then the elements
\begin{equation}
A_a=K_a{}^{ij}\left(\pi_{ij}\ot\id\right)\,A\,\in\uqg
\end{equation}
span a weak quantum Lie algebra $\wqlie=\mbox{span}_{\ch}
\{A_a\}$.
\end{proposition}

\begin{proof}
The expression $A=\hhh^{-1}\left(R^T R-1\right)$ is well defined
because $R=1 \modh$. It follows from the defining property 
$R\,\Delta(x)=\Delta^T(x)\,R\ \forall x\in\uqg$
of the R-matrix that
$A\,\Delta(x)=\Delta(x)\,A,~~\forall x\in\uqg$. It is then clear
from \lemref{lem} that the $A_a$ span an ad-submodule
of $\uqg$. It follows from the definition of the Clebsch-Gordan
coefficients $K_a{}^{ij}$ that this ad-submodule is either
isomorphic to the adjoint representation or zero.
$R$ satisfies $R=1+\hhh\,r+{\cal O}(\hhh^2)$ where $r\in\lie\ot\lie$
is the classical $r$-matrix. Thus $A=r+r^T \modh\in\lie\ot\lie$
and $A_a\modh\in\lie$. It follows that $\mbox{span}_{\ch}
\{A_a\}=\lie\modh$.
\end{proof}

Using the fact, established in \cite{qlie}, that given a weak
quantum Lie algebra $\wqlie$ one can always construct a 
true quantum Lie algebra $\lqlie$, we arrive at the announced
existence result.

\begin{theorem}
For any finite-dimensional simple complex Lie algebra $\lie$
there exists at least one quantum Lie algebra $\lqlie$ inside $\uqg$.
\end{theorem}

\noindent{\bf Thanks:}
We thank Andrew Pressley, Vyjayanthi Chari, Manfred Scheunert and Chris 
Gardner for discussions and helpful comments.


\begin{thebibliography}{99}

\bibitem{toda} H.W. Braden, E. Corrigan, P.E. Dorey, R. Sasaki, 
 {\it Affine Toda field theory and exact S-matrices}, 
 Nucl. Phys. {\bf B338} (1990) 689; 
G.W. Delius, M.T. Grisaru, D. Zanon, {\it Exact S-matrices
for nonsimply-laced affine Toda theories}, 
hep-th/9201067, Nucl. Phys. {\bf B382} (1992) 365.

\bibitem {Dri85}        V.G. Drinfel'd,  {\it Hopf algebras and the quantum
   Yang-Baxter equation}, Sov. Math. Dokl. {\bf 32} (1985) 254.

\bibitem {Jim85}     M. Jimbo,   {\it A q-Difference Analogue of U(g) and
  the Yang-Baxter Equation}, Lett. Math. Phys. {\bf 10} (1985) 63.

\bibitem{qlie} G.W. Delius, A. H\"uffmann, {\it On Quantum Lie Algebras
 and Quantum Root Systems}, q-alg/9506017, J. Phys. A {\bf 29} (1996) 1703.

\bibitem{Wor89} S.L. Woronowicz, {\it Differential Calculus on Compact
Matrix Pseudo\-groups (Quantum Groups)}, Comm. Math. Phys. {\bf 122}
(1989) 125.

\bibitem {difcalc} P. Aschieri, L. Castellani, {\it An introduction to
noncommutative differential geometry on quantum groups},
Int. J. Mod. Phys. {\bf A8} (1993) 1667;
\vspace{1mm}\hfill\newline
D. Bernard,      {\it Quantum Lie Algebras and
Differential Calculus on Quantum Groups}, Prog. Theo. Phys. Suppl. {\bf 102}
(1990) 49;
\vspace{1mm}\hfill\newline
B. Jurco,       {\it Differential Calculus on Quantized Simple
Lie Groups}, Lett. Math. Phys. {\bf 22} (1991) 177;
\vspace{1mm}\hfill\newline
P. Schupp, P. Watts, B. Zumino,        {\it Bicovariant
Quantum Algebras and Quantum Lie Algebras}, Commun. Math. Phys. {\bf 157}
(1993) 305;
\vspace{1mm}\hfill\newline
P. Schupp, {\it Quantum Groups, Non-Commutative Differential
Geometry and Applications}, hep-th/9312075 (1993).

\bibitem {Sud95} A. Sudbery, V. Lyubashenko, {\it Quantum Lie Algebras of Type} $A_n$, 
q-alg/9510004.

\bibitem {Sch} K. Schm\"udgen and A. Sch\"uler, {\it Left-covariant 
Differential Calculi on $SL_q(N)$}, Leipzig preprint.


\bibitem {diffcalc} G.W. Delius, {\it The Problem of Differential 
Calculus on Quantum Groups}, q-alg/9608019.

\bibitem {Cha94}        V. Chari, A. Pressley,   {\it A Guide to Quantum
 Groups}, Cambridge University Press {\bf } (1994).

\bibitem {Cur62} C.W. Curtis, I. Reiner, {\it Representation theory
of finite groups and associative algebras}, Interscience Publishers (1962).

\bibitem{Swe69} M.E. Sweedler, {\it Hopf algebras}, Benjamin, New York (1969).

\bibitem{Dri90} V.G. Drinfel'd, {\it On almost cocommutative Hopf algebras},
Leningrad Math. J. {\bf 1} (1990) 321.

\bibitem{Din94} J. Ding, I.B. Frenkel, {\it Spinor and Osciallator
Representations of Quantum Groups}, Prog. Math. {\bf 123} (1994) 127.

\bibitem {qgln} G.W. Delius, M.D. Gould, A. H\"{u}ffmann, Y.-Z. Zhang,
 {\it Quantum Lie algebras associated to $U_q(gl_n)$ and $U_q(sl_n)$},
 q-alg/9508013, J. Phys. A 29 (1996) 5611.

\bibitem {Res87} N.Yu. Reshetikhin, {\it Quantized Universal Enveloping 
Algebras, the Yang-Baxter Equation and Invariants of Links. I.} 
LOMI-preprint E-4-87.


\end{thebibliography}
\end{document}